\documentclass{amsart}

\newtheorem{theorem}{Theorem}[section]
\newtheorem{lemma}[theorem]{Lemma}
\newtheorem{proposition}[theorem]{Proposition}
\newtheorem{corollary}[theorem]{Corollary}
\newtheorem{conjecture}[theorem]{Conjecture}

\theoremstyle{definition}
\newtheorem{definition}[theorem]{Definition}

\theoremstyle{remark}
\newtheorem{remark}[theorem]{Remark}

\numberwithin{equation}{section}

\def\P#1{{\mathbb P}^#1}

\newcommand{\s}{\mathcal}
\newcommand{\sL}{{\s L}}
\newcommand{\M}{{\s M}}
\newcommand{\N}{{\s N}}
\newcommand{\sO}{{\s O}}

\newcommand{\propref}[1]{Proposition~\ref{#1}}
\newcommand{\thmref}[1]{Theorem~\ref{#1}}
\newcommand{\lemref}[1]{Lemma~\ref{#1}}
\newcommand{\corref}[1]{Corollary~\ref{#1}}

\newcommand{\conjref}[1]{Conjecture~\ref{#1}}

\begin{document}

\title[Non-Defectivity of Grasmannians of planes]{Non-Defectivity of Grassmannians of planes}
\date{June 26, 2008}
\author{Hirotachi Abo}
\email{abo@uidaho.edu}
\address{Department of Mathematics, University of Idaho, Moscow, ID 83844, USA}
\author{Giorgio Ottaviani}
\email{ottavian@math.unifi.it}
\address{Dipartimento di Matematica ``Ulisse Dini", Universit\`a degli Studi di Firenze, 
Viale Morgagni 67/A, 50134 Firenze, Italy}
\author{Chris Peterson}
\email{peterson@math.colostate.edu}
\address{Department of Mathematics, Colorado State University, Fort Collins, CO 80523-1874, USA}
\date{}

\begin{abstract}
Let $Gr(k,n)$ be the Pl\"ucker embedding of the Grassmann variety of projective $k$-planes in $\P n$. For a projective variety $X$, let  $\sigma_s(X)$ denote the variety of its $s-1$ secant planes. More precisely, $\sigma_s(X)$ denotes the Zariski closure of the union of linear spans of $s$-tuples of points lying on $X$. We exhibit two functions
$s_0(n)\le s_1(n)$ such that $\sigma_s(Gr(2,n))$ has the expected dimension whenever $n\geq 9$ and either
$s\le s_0(n)$ or $s_1(n)\le s$. Both $s_0(n)$ and $s_1(n)$ are asymptotic to $\frac{n^2}{18}$. This yields, asymptotically, the typical rank of an element of $\wedge^{3}\hskip 1pt {\mathbb C}^{n+1}$.
Finally, we classify all defective $\sigma_s(Gr(k,n))$ for $s\le 6$ and provide geometric arguments underlying each defective case.
\end{abstract}

\maketitle

\section{Introduction}

Let  $X \subset \P N$ be a non-degenerate projective variety.  The {\it $s$-secant variety} $\sigma_s(X)$ is defined to be the Zariski closure of the union of linear spans of $s$-tuples of points lying on $X$ (see \cite{Z}). 
Note that with this notation, $\sigma_2(X)$ is the usual variety of secant lines of $X$. There is a smallest $s$ such that $\sigma_s(X)=\P N$ leading to a natural filtration: 
\begin{eqnarray*}
\label{seq:filtration}
X=\sigma_1(X)\subset\sigma_2(X)\subset\sigma_3(X)\subset\cdots \subset \sigma_s(X) =\P N.
\end{eqnarray*}

Let $Gr(k,n)$ denote the Grassmannian of projective $k$-planes in $\P n$. For the purposes of this paper, we will assume that $Gr(k,n)$ is embedded through the Pl\"ucker map in $\P N$ with $N={{n+1}\choose{k+1}}-1$. We can identify points in $\P N$ with general skew-symmetric tensors and points on $Gr(k,n)$ as decomposable skew-symmetric tensors.  An element $\omega\in \wedge^{k+1}\mathbb{C}^{n+1}$ has {\it rank r} if it can be written as a linear combination of 
$r$ decomposable skew-symmetric tensors (but not fewer). In other words, $\omega=\sum_{t=1}^r v_{1,t}\wedge\ldots\wedge v_{k+1,t}$ with $v_{i,j}\in{\mathbb C}^{n+1}$. The higher secant variety $\sigma_s(Gr(k,n))$  can be viewed as a compactification of the ``parameter space" for skew-symmetric tensors of rank less than or equal to $s$. An interesting problem related to the rank of skew-symmetric tensors is to find the least integer $\underline R(k,n)$ such that a generic skew-symmetric tensor has rank less than or equal to $\underline R(k,n)$. The integer $\underline R(k,n)$ is called the {\it typical rank} of  $\wedge^{k+1}\mathbb{C}^{n+1}$ (also called  the {\it essential rank} in \cite{E}). Note that the filtration of skew-symmetric tensors by their ranks leads naturally to an identification of $\underline R(k,n)$ as the least integer $s$ such that $\sigma_s(Gr(k,n))=\P N$. See \cite{LM} for a recent survey on the subject and its applications.

If  $k=1$ then $X$ is a Grassmannian of lines and $\sigma_s(X)$ corresponds to the locus of skew-symmetric morphisms of rank less than or equal to $s$. It is well known that
a skew-symmetric morphism, corresponding to a skew-symmetric matrix of rank $2s$, can be written as the sum of $s$ decomposable skew-symmetric tensors (but not fewer).  In particular, we have  $\underline R(1,n)=\lceil\frac{n+1}{2}\rceil$. Thus we may assume
that $k\ge 2$.

It is straightforward to show that $$\frac{{{n+1}\choose{k+1}}}{(k+1)(n-k)+1}\le\underline R(k,n)$$
(see the next section). In particular, we have 
$$\frac{n^2}{18}+O(n)\le\underline R(2,n).$$

On the other hand, Ehrenborg found in \cite{E}, Corollary 7.9, the  upper bound
$$\underline R(2,n)\le \frac{(n^2+3)}{12}+1$$
by using results on Steiner triple systems.
One of the main goals of this paper is to give a sharp asymptotic bound for $\underline R(2,n)$. To be more precise, we will prove the following theorem: 
\begin{theorem}\label{main} If $\underline R(2,n)$ denotes the rank of a generic skew symmetric tensor $\omega\in \wedge^{3}\mathbb{C}^{n+1}$
then $\underline R(2,n)\sim \frac{n^2}{18}.$
\end{theorem}
Theorem \ref{main} is obtained as a consequence of a more precise (but more technical) theorem on the dimension of $\sigma_s(Gr(k,n))$. 
To state this theorem, we will need to review several known facts from the literature.

First of all, the following inequality is easy to establish (see the next section for a geometric interpretation): 
\[
\dim\sigma_s(Gr(k,n))\leq \min\left\{s[(k+1)(n-k)+1]-1, \ \frac{{{n+1}\choose{k+1}}}{(k+1)(n-k)+1}\right\}. 
\]
We say that $\sigma_s(Gr(k,n))$ has the {\it expected dimension} if equality holds.
If there exists a $s$ for which $\sigma_s(Gr(k,n))$ does not have the expected dimension then $Gr(k,n)$ is  said to be {\it defective}. As we have seen, Grassmannians of lines are nearly always defective. In contrast, 
there are very few cases where $Gr(k,n)$ is known to be defective when $k\geq 2$. 

In 1916, C. Segre \cite{Se} 
proved that $\sigma_2(Gr(2,5))$
has the expected dimension which established that  $\underline R(2,5)=2$.
On the other hand, in 1931, Schouten \cite{Sch} showed that $Gr(2, 6)$ is defective. 
Indeed he proved that $\sigma_3 (Gr( 2, 6))$ is a hypersurface as opposed to filling the ambient space. This result established that $\underline R(2,6)=4$.
It is well known that the degree of Schouten's hypersurface is seven (\cite{La}). In Section 5
we analyze this case in more detail. In particular, we find an explicit description of this degree seven invariant by relating its cube to the
determinant of a $21\times 21$ symmetric matrix.

\begin{theorem}\label{t1}
Let $\omega\in\wedge^3{\mathbb C}^7$. Consider the contraction operator $\phi_{\omega}\colon\wedge^2{\mathbb C}^7\to\wedge^5{\mathbb C}^7. $
The equation of $\sigma_3(Gr(2,6))$ is given by an $SL(7)$-invariant polynomial $P_7$
of degree seven such that
$$\det(\phi_{\omega})=2\left[P_7(\omega)\right]^3.$$
\end{theorem}

In 2002, Catalisano, Geramita and Gimigliano (with the help of Catalano-Johnson who had some unpublished results on this subject) \cite{CGG1} showed that $Gr(3, 7)$ and $Gr(2, 8)$ are defective. Due to the isomorphism $Gr(k,n)\cong Gr(n-k-1,n)$ (for instance, $Gr(2,8)\cong Gr(5,8)$), we only consider Grassmannians, $Gr(k,n)$, for which $k\leq \frac{n-1}{2}$. Based on a mixture of theory and computational experiments there is a body of evidence suggesting that all defective Grassmannians have been found. As a result, we believe in the following conjecture proposed in \cite{BDG} (conj. 4.1):
\begin{conjecture}[Baur-Draisma-de Graaf]  \label{BDdG}
Let $k\ge 2$. Then $\sigma_s(Gr(k,n))$ has the expected dimension except for the following cases:
 $$\begin{array}{c|c|c|c}
 &&\textrm{actual codimension}&\textrm{expected codimension}\\
 \hline
 (1)&\sigma_3(Gr(2,6))&1&0\\
 \hline (2)&\sigma_3(Gr(3,7))&20&19\\
  \hline (2')&\sigma_4(Gr(3,7))&6&2\\
  \hline (3)&\sigma_4(Gr(2,8))&10&8\\
 \end{array}$$
\end{conjecture}
If the conjecture is true, then $\sigma_3(Gr(2,6))$ is the only secant variety which both does not have the expected dimension and is a hypersurface (with $k\geq 2$). The invariant computed in \thmref{t1} defines this hypersurface as its zero-locus.
Computational evidence was given by McGillivray who performed a Montecarlo technique to check that the conjecture is true for $n\le 14$ in \cite{McG}. This result was extended to $n=15$ in \cite{BDG}. One of the goals of this paper is to provide further evidence in support of this conjecture.  
As a step in this direction, in Section 3 we classify  Grassmann varieties with defective $s$-secant varieties for small $s$. More precisely, we prove the following theorem:

\begin{theorem}\label{t2} Except for the cases listed in Conjecture \ref{BDdG},  $\sigma_s(Gr(k,n))$ has the expected dimension whenever $k\geq 2$ and $s\le 6$.
\end{theorem}

Let $A(n,6,w)$ be the cardinality of the largest binary code
of length $n$,  constant weight $w$, and distance 6 (see Section $3$ for more details).
It is not hard to show that if $s\leq A(n+1,6,k+1)$ then
$\sigma_s(Gr(k,n))$ has the expected dimension.
 For small values of $n$ and $k$ this is a useful result via the monomial approach, like in \cite{E} (indeed we use it
in Section 3). However,
for $n\gg 0$ the value of $A(n+1,6,k+1)$ is typically
smaller than $$\max\left\{ s \ \left|\ \sigma_s(Gr(k,n))\ \textrm{\ has the expected dimension and does not fill}\ \P N\right. \right\}.$$
For example $A(10,6,4)=5$ by \cite{S}, while by \thmref{t2} we see that
$\sigma_6(Gr(3,9))$ has the expected dimension and does not fill the ambient space.

We can now state a slightly technical theorem which implies \thmref{main}.
We show that there are two functions
$s_1(n)\le s_2(n)$ such that  $\sigma_s(Gr(2,n))$ has the expected dimension whenever
either $s\le s_1(n)$ or $s\ge s_2(n)$. The precise statement is the following: 

\begin{theorem}\label{asymp}
Let $n\ge 9$. Let 
\[
s_1(n)=\left\lfloor\frac{n^2}{18}-\frac{20n}{27}+\frac{287}{81}\right\rfloor+
\left\lfloor\frac{6n-13}{9}\right\rfloor
\]
and let 
\[
 s_2(n)=\left\lceil\frac{n^2}{18}-\frac{11n}{27}+\frac{44}{81}\right\rceil+
\left\lceil\frac{6n-13}{9}\right\rceil. 
\]
Then $\sigma_s(Gr(2,n))$ has the expected dimension whenever $s\le s_1(n)$
and whenever $s\ge s_2(n)$ (in this second case it fills the ambient space).
\end{theorem}

Ideally, we would like the functions to satisfy
$s_1+1 \geq s_2$ modulo a finite list of exceptions. While the theorem does not reach this result, it does have a relatively small value for $s_2-s_1$. Such a result is reminiscent of one
obtained in \cite{CGG2} for $X=\P 1\times\ldots\times\P 1$  where they produced functions 
$s_1(t),s_2(t)$ (with $t$ denoting the number of factors of $\P 1$) such that $s_2-s_1\le 1$. This was extended in \cite{AOP1}
to $X=\P n\times\ldots\times\P n$, where the functions satisfy
 $s_2-s_1\le n$. In the recent \cite{CGG3} the final result for $X=\P 1\times\ldots\times\P 1$ 
 has been found.
Note that in \thmref{asymp}, $s_1(n)\sim s_2(n)\sim\frac{n^2}{18}$ (the sharp asymptotical value) 
and that
\thmref{main} follows.  

The proof of \thmref{asymp} is in Section 4. Our approach relies on a specialization technique to place a certain number of points on subgrassmannians determined by codimension six linear subspaces (see the remark after
\propref{propa}). 
This technique was inspired by \cite{BO}, where the case $X=(\P n,\sO(3))$ was treated
with a similar specialization determined by codimension three linear subspaces.

The present technique can be extended
to higher values of $k$ (see \cite{AOP2}), but at the price of a much more complicated inductive procedure involving joins,
and with a massive use of the computer to check initial cases. As a consequence, we
decided to treat this case in a separate paper.

In Section 5 we close the paper with a geometric explanation for each of the known defective cases appearing in the list of Conjecture \ref{BDdG}.
\section{Notation and definitions}
We begin this section by recalling the precise definition of a higher secant variety.
 If $Q_1,\dots, Q_s$ are points in $\P m$ then we let $\langle Q_1,\dots, Q_s \rangle$ denote their linear span. If $X\subseteq\P m$ is
a projective variety then the {\it $s$-secant variety} $\sigma_s(X)$ is defined to be the Zariski closure of the union of
 the linear span of $s$-tuples of points $(Q_1,\dots, Q_s)$ where $Q_1,\dots Q_s\in X$. 
 In other words $$\sigma_s(X)=\overline{\bigcup_{Q_1, \dots, Q_s\in X} \langle Q_1,\dots,Q_s\rangle}.$$
 If $X\subset \P m$ is a non-degenerate variety of dimension $d$ then a standard dimension count leads to an upper bound on the dimension of $\sigma_s(X)$. In particular, the dimension of $\sigma_s(X)$ can never exceed $\min\{s(d+1)-1,m\}$. Determining when the secant variety of a Grassmann variety reaches this upper bound is one of the main goals of this paper. We summarize some of the terminology that will be used to develop the main theorems in the following:
\begin{definition}
Let $X$ be a non-degenerate $d$-dimensional variety in $\P m$.
\begin{itemize}
\item[(1)] If $\dim \sigma_s(X)=\min\{s(d+1)-1,m\}$ then $\sigma_s(X)$ is said to have the {\it expected dimension}.
\item[(2)] If $\dim\sigma_s(X)<\min\{s(d+1)-1,m\}$ then $X$ is said to have a {\it defective $s$-secant variety}.
\item[(3)] If there exists an $s$ such that $\dim\sigma_s(X)<\min\{s(d+1)-1,m\}$ then $X$ is said to be {\it defective}.
\item[(4)] The smallest integer $s$ such that $\sigma_s(X)$ fills the ambient space
is called the {\it typical rank} and is denoted by $\underline R(X)$.
\end{itemize}
\end{definition}

The main tool that will be used to compute the dimension of $\sigma_s(X)$ is the following celebrated theorem of Terracini: 
\begin{theorem}[Terracini's Lemma \cite{Z}] 
Let $P_1,\dots, P_k$ be points in $X$ and let $z$ be a general point in $\langle P_1,\ldots ,P_k \rangle$. Then the tangent space to $\sigma_s(X)$ at $z$ is given by
   $$T_z\sigma_s(X)=\langle T_{P_1}X,\ldots ,T_{P_k}X \rangle$$ where $T_{P_i}X$ denotes the tangent space to $X$ at $P_i$.
\end{theorem}

Let $\mathbb{K}$ be a field with $\mathrm{char}({\mathbb K})=0$ and 
let $V$ be an $(n+1)$-dimensional 
vector space over $\mathbb{K}$.   
We denote by $Gr(k,n)$ the Grassmannian 
of $(k+1)$-dimensional subspaces of $V$ and by $C(Gr(k,n))$ 
the affine cone over $Gr(k,n)$. 

We denote by  
$\lfloor x\rfloor$ the greatest integer less than or equal to $x$ and by $\lceil x\rceil$ the smallest integer greater than or equal to $x$.
\section{Classification of Grassmannians with defective $s$-secant varieties, $s\le  6$}
In this section we classify Grassmannians with defective $s$-secant varieties when $s\le 6$. The main tools are a combination of computations on the computer and the monomial technique. We use ideas from coding theory to strengthen the monomial approach. Throughout this section,  $\mathbb{K}$ is an infinite field with $\mathrm{char}\  (\mathbb{K})\not=2$.

 We begin with two well known lemmas.
\begin{lemma}[Tangent space at a point of a Grassmannian]
\label{tangent}
 Let $p=v_0\wedge\ldots\wedge v_k$
be a point of
$C(Gr(k,n))$ where $v_i\in V=\mathbb{K}^{n+1}$. The tangent space to $C(Gr(k,n))$ at $p$ is
\[
T_p(k,n)=\sum_{i=1}^k v_1 \wedge \cdots \wedge v_{i-1} \wedge V \wedge v_{i+1} \wedge \cdots \wedge v_k.
\]   
\end{lemma}

Let ${\mathcal B}=\{v_0,\dots ,v_n\}$ be a basis of $V=\mathbb{K}^{n+1}$. The ambient space of $C(Gr(k,n ))$ in its
 Pl\"ucker embedding is $\wedge^{k+1}V$. A basis of the ambient space is determined by the $(k+1)$-element subsets of ${\mathcal B}$.
 From \lemref{tangent} we see that a basis of the tangent space to $C(Gr(k,n))$ at
 $v_0\wedge\ldots\wedge v_k$ is determined by the $(k+1)$-element subsets of ${\mathcal B}$ which intersect
  $\{v_0,\dots, v_k\}$ in a set of at least $k$ elements. Thus, we have the following Lemma:

\begin{lemma}[Monomial Lemma]
\label{monomial}  Let $\{v_0,\dots, v_n\}$ be a basis for $V$ and let $k\geq 2$. 
\begin{itemize}
\item[(1)] Let ${\mathbf A}=\{a_0,\dots ,a_k\}$ be a subset of $\{v_0,\dots ,v_n\}$ and let $T_{\mathbf A}(k,n)$ denote the tangent space to $C(Gr(k,n))$ at $a_0\wedge\ldots\wedge a_k$. A basis of $T_{\mathbf A}(k,n)$ is given by vectors corresponding to $(k+1)$-element subsets of $\{v_0,\dots, v_n\}$ which intersect ${\mathbf A}$ in at least $k$ elements.
\item[(2)] Let ${\mathbf A}_1,\dots, {\mathbf A}_t$ be $(k+1)$-element subsets of $v_0,\dots, v_n$. If $|{\mathbf A}_i\cap {\mathbf A}_j|\leq k-2$ whenever $i\neq j$ then $T_{\mathbf A_1}(k,n),\dots, T_{\mathbf A_t}(k,n)$ are linearly independent.
\end{itemize}
\end{lemma}

By upper semicontinuity, if there exist smooth points $Q_1, \dots, Q_s \in Gr(k,n)$ such that the tangent spaces $T_{Q_1}(k,n), \dots, T_{Q_s}(k,n)$ are linearly independent (or else span the ambient space) then $\sigma_s(Gr(k,n))$ has the expected dimension by Terracini's Lemma. Thus, to show that $X$ does not have a defective $s$-secant variety, it is enough to find $s$ smooth points on $X$ such that the tangent spaces at these points are either linearly independent or else span the ambient space. 

The following theorem extends Theorem 2.1 of \cite{CGG1}.
\begin{theorem}\label{tre}
If $3(s-1)\le n-k$ and if $k\ge 2$
then $$\dim\sigma_s(Gr(k,n))=s(k+1)(n-k)+(p-1).$$
\end{theorem}
\begin{proof}
 Let $\P {n}={\mathbb P}(V)$. Fix a basis $\{v_0,\ldots,v_n\}$ for $V$.
Let $\mathbf{A}_1,\ldots, \mathbf{A}_s$ be points in $Gr(k,n )$ defined by $\mathbf{A}_i=\langle v_{3(i-1)},\ldots, v_{3(i-1)+k}\rangle$.
By \lemref{monomial}, it follows that the tangent spaces at $\mathbf{A}_1,\ldots , \mathbf{A}_s$ are linearly independent, hence by Terracini's lemma we are done.
\end{proof}

\begin{remark} By \lemref{monomial} and Terracini's lemma, we can show that $\sigma_s(Gr(k,n))$  has the expected dimension if we can show that there exist $s$ distinct $(k+1)$-element subsets, ${\mathbf A}_1,\dots, {\mathbf A}_s$, of an $(n+1)$-element set such that whenever $i\neq j$ we have $|{\mathbf A}_i\cap {\mathbf A}_j|\leq k-2$. To each $(k+1)$-element subset, we can associate a weight $k+1$ binary vector of length $n+1$ via the characteristic function. Our conditions require that the Hamming distance between any pair of distinct vectors is at least 6 (this was observed also in \cite{BDG}). Let $A(n,6,w)$ denote the cardinality of the largest set of length $n$, weight $w$ vectors that satisfy this condition on the Hamming distance. Lower bounds on $A(n,6,w)$ have been computed as part of a search for good {\it constant weight binary codes}. Tables of such bounds can be used to prove that certain Grassmann varieties are not $s$-defective via monomial methods. We found the table \cite{S} and the paper \cite{GSl} particularly useful.
\end{remark}

\begin{theorem}[Graham and Sloane \cite{GSl}] \label{graham-sloane}  Let $A(n,6,w)$ denote the maximum number of codewords in any binary code of length $n$, constant weight $w$ and Hamming distance 6.
\begin{itemize}
\item[(a)] Let $q$ be the smallest prime power such that  $q\ge n$, then $A(n,6,w)\ge \frac{1}{q^2}$ ${n}\choose{w}$.
\item[(b)] Let $q$ be the smallest prime power such that  $q+1\ge n$, then $A(n,6,w)\ge \frac{q-1}{q^3-1}$ ${n}\choose{w}$.
\item[(c)] $A(n,6,w)\ge$ ${n}\choose{w}$ $/ $ $( {1+w(n-w)+ {{w}\choose{2}}{{n-w}\choose{2}}})$.
\end{itemize}
\end{theorem}

Using results of \cite{McG}, \cite{S} and applying \thmref{tre} and \thmref{graham-sloane},  we obtain the following corollary (which proves \thmref{t2}): 
\begin{corollary}\label{p6} If $k\geq 2$ then $\sigma_s(Gr(k,n))$ satisfies the following:
\begin{itemize}
\item[(i)] $\sigma_3(Gr(k,n))$ has the expected dimension
except for $(k,n)=(2,6), (3,7)$
\item[(ii)] $\sigma_4(Gr(k,n))$ has the expected dimension
except for $(k,n)=(2,8), (3,7)$
\item[(iii)]  $\sigma_s(Gr(k,n))$ always has the expected dimension for $s=2,5,6$.
\end{itemize}
\end{corollary}
\begin{proof} The case $s=2$ was established in Corollary 2.2 of \cite{CGG1}.

For the case $s=3$,  we apply  \thmref{tre} when the inequality $n-k\ge 6$ is satisfied. The remaining cases have $n-k\le 5$
with $n\le 9$ (due to $k\le \frac{n-1}{2}$). These have been checked in \cite{McG}.

 For the case $s=4$, we apply \thmref{tre} when the inequality $n-k\ge 9$ is satisfied. 
 The remaining cases have $n-k\le 8$
and $n\le 15$. For $n\le 14$, they have been checked in \cite{McG}. The only remaining case is $(k,n)=(7,15)$ which
follows from \cite{S}.

For the case $s=5$, we apply \thmref{tre} when the inequality $n-k\ge 12$ is satisfied.
 The remaining cases have $n-k\le 11$ and $n\le 21$. For $n\le 14$ they have been checked in \cite{McG}.
 The only remaining cases have $k\ge 4$ and $21\ge n\ge 15$. These all follow from \cite{S}.

For the case $s=6$, we apply \thmref{tre} when the inequality $n-k\ge 15$ is satisfied. 
 The remaining cases have $n-k\le 14$ and $n\le 27$.   The case $(k,n)=(2,16)$ can be checked by the
computer exactly as in \cite{McG}. The remaining cases have $k\ge 3$ and $15\le n \le 27$. These all follow from \cite{S}.
\end{proof}

\begin{remark} 
By using \thmref{tre}, \thmref{graham-sloane},
the table \cite{S}, and the algorithm in \cite{McG}, it is expected that one can extend \corref{p6} to larger values of $s$.
\end{remark}

\section{The inductive step, from $n-6$ to $n$}
In this section we develop a collection of tools that lead to an inductive proof of \thmref{asymp}. Throughout the section, 
we denote by $V$ an $(n+1)$-dimensional vector space over an infinite field $\mathbb{K}$ and we denote
by $V^*$ its dual space. 
\begin{proposition}\label{propa}
Let $X=Gr(2,n)$ with $n\ge 17$. Let $V=\mathbb{K}^{n+1}$, and let $L$, $M$ and $N$ be general codimension six subspaces of $V$. 
Let $\sL$, $\M$ and $\N$ be the Grassmann varieties of $3$-planes in $L$, $M$ and $N$ respectively. Let  $p_1,\dots ,p_4$ be $4$ general points on $\sL$,  $q_1,\dots ,q_4$ be $4$ general points on $\M$, and  $r_1,\dots ,r_4$ be $4$ general points on $\N$. Then there are no hyperplanes in $\P{{}}(\wedge^3 V)$ which contain both
$\sL\cup\M\cup\N$ and the tangent spaces to $X$ at the $12$ points $\{p_i,q_i,r_i\}_{1 \leq i \leq  4}$.
\end{proposition}
\begin{proof} Let $\{e_0, \dots, e_n\}$ be a basis for $V$ 
and let $\{x_0, \dots ,x_n\}$ be its dual basis. Without loss of generality, we may assume that $L=\{x_i=0, i=0,\ldots,5\}$,
$M=\{x_i=0, i=6,\ldots,11\}$ and 
$N=\{x_i=0, i=12,\ldots,17\}$.
The hyperplanes in $\P{{}}(\wedge^3 V)$ which contain
$\sL\cup\M\cup\N$ span a space of dimension $6^3=216$ with basis
$\{x_0\wedge x_6\wedge x_{12},\ldots , x_5\wedge x_{11}\wedge x_{17}\}$.

We remark that the codimension of $\sL$ (resp. $\M$, $\N$) in $X$ is $18$.
Furthermore, $12 \cdot [3(n+1-3)-3(n-5-3)]=12\cdot 18=216$. 
Thus it is enough to prove that 
the tangent spaces to $X$ at the points $p_1,p_2,p_3,p_4,q_1,q_2,q_3,q_4,r_1,r_2,r_3,r_4$ modulo 
$\langle \P{{}}(L), \P {{}}(M), \P{{}}(N) \rangle$ form a $216$-dimensional vector space.  To prove this, it is enough to do the computation in $\P {{17}}$ which we achieve
 through the Macaulay2 script given below:
\begin{footnotesize}
\begin{verbatim}
i1 : randomIdeal = method();
i2 : randomIdeal(Ideal) := Ideal => I -> (
          R := ring I; 
          ideal(gens I*random(source gens I,R^{3:-1}))
          ); 
i3 : randomIdeal' = method(); 
i4 : randomIdeal'(Ideal) := Ideal => I -> (
          R := ring I;
          J := I^2;
          m := gens J; 
          ideal(m*map(source m,,basis(3,J)))
          ); 
i5 : E = ZZ/32003[e_0..e_17, SkewCommutative=>true];
i6 : L = ideal(e_6..e_17);
o6 : Ideal of E
i7 : M = ideal(e_0..e_5,e_12..e_17);
o7 : Ideal of E
i8 : N = ideal(e_0..e_11);
o8 : Ideal of E
i9 : lt = {L,M,N};
i10 : h = new MutableList; 
i11 : scan((#lt),i->h#i=(lt#i)^3); 	
i12 : J = trim(sum(toList h));
o12 : Ideal of E
i13 : I = ideal(0_E); 
o13 : Ideal of E
i14 : for i from 1 to 4 do (
          h' = new MutableList;
           scan((#lt),i->h'#i=randomIdeal'(randomIdeal(lt#i)));
           I = I+sum(toList h');
           );
i15 : rank source gens trim (I+J) == rank source gens trim (ideal vars E)^3 
o15 = true
\end{verbatim}
\end{footnotesize}
Note that we work in characteristic $p=32003$.   Our goal is to check that a certain integer matrix has maximal rank. The Macaulay2 script determines that the matrix has maximal rank modulo $p$. 
The result in characteristic zero follows from 
the  openess of the  maximal rank condition. 
\end{proof}

\begin{remark} It would be natural for the reader to ask why we choose subspaces of codimension six.
 The linear system of hyperplanes in $\P{{}}(\wedge^3 V)$ which contain
the union of three subspaces of codimension $p$ has dimension $p^3$ when $n$ is sufficiently large.
 Any tangent space supported on this union imposes
$3p$ conditions. In order to have a number of points which impose independent conditions
on the linear system and which make it empty, we need the condition
that $3p$ divides $p^3$. Therefore, a necessary condition is that $p$ is a multiple of $3$. Unfortunately,
the case $p=3$ does not work. Three general points supported on three codimension three
subspaces do not impose independent conditions, as can be
checked by an analogous Macaulay2 script. Hence we turned to the next case, $p=6$.
\end{remark}

Let $f(n)={{n+1}\choose 3}$. We compute the finite difference $f(n)-2f(n-6)+f(n-12)=36(n-6)$.
In particular the system of hyperplanes in $\P{{}}(\wedge^3 V)$ which contain
$\sL\cup\M$ has dimension $36(n-6)$ (by the Grassmann formula) for $n\ge 11$.
In order to fit with \propref{propa}, we remark that
$f(n)-3f(n-6)+3f(n-12)-f(n-18)=216$.

We want to keep four points outside $\sL\cup\M$. Note
that the tangent space, at each of the four points, imposes $3n-5$ conditions
and that
$$\frac{36(n-6)-4\cdot(3n-5)}{36}=\frac{6n-49}{9}.$$ This leads to the following proposition:

\begin{proposition}\label{propb}
Let $X=Gr(2,n)$ with $n\ge 11$. Let $V=\mathbb{K}^{n+1}$, and let $L$ and $M$ be general codimension six subspaces of $V$. 
Let $\sL$ (resp. $\M$) be the Grassmann variety of $3$-planes in $L$ (resp. $M$). Then

\begin{itemize}
\item[(i)] The system of hyperplanes in $\P{{}}(\wedge^3 V)$ which contain
$\sL\cup\M$ and which contain the tangent spaces at $\lfloor\frac{6n-49}{9}\rfloor$
general points on $\sL$, at $\lfloor\frac{6n-49}{9}\rfloor$
general points on $\M$,
and at $4$ general points on $X$
has the expected dimension $36(n-6)-36\lfloor\frac{6n-49}{9}\rfloor-4(3n-5)$,
which is
$$\left\{\begin{array}{ccc}20&\textrm{if}&n=0\ (\textrm{mod\ }3)\\
8&\textrm{if}&n=1\ (\textrm{mod\ }3)\\
32&\textrm{if}&n=2\ (\textrm{mod\ }3)\\
\end{array}\right.$$

\item[(ii)]  There are no hyperplanes in $\P{{}}(\wedge^3 V)$ which contain
$\sL\cup\M$ and which contain the tangent spaces at $\lceil\frac{6n-49}{9}\rceil$
general points on $\sL$, at $\lceil\frac{6n-49}{9}\rceil$
general points on $\M$ 
and at $4$ general points on $X$.
\end{itemize}
\end{proposition}
\begin{proof}
We will let $\{p_i\}$ denote the $\lfloor\frac{6n-49}{9}\rfloor$ (resp. $\lceil\frac{6n-49}{9}\rceil$) general points on $\sL$, $\{q_i\}$ denotes the $\lfloor\frac{6n-49}{9}\rfloor$ (resp. $\lceil\frac{6n-49}{9}\rceil$) general points on $\M$ and $\{r_i\}$ denotes the four general points on $X$.
The proof is by a 6-step induction from $n-6$ to $n$. The initial six cases
$n=11, 12, 13, 14, 15, 16$ can be checked directly as follows:
\begin{footnotesize}
\begin{verbatim}
i16 : secondStep = method()     
o16 = secondStep
o16 : MethodFunction
i17 : secondStep(ZZ) :=  n -> (
           s := floor((6*n-49)/9);
           t := binomial(n+1,3)-(36*(n-6)-36*s-4*(3*n-5)); 
           E := ZZ/32003[e_0..e_n, SkewCommutative=>true];
           lt := {ideal(e_6..e_n),ideal(e_0..e_(n-6))};
           J := sum(apply(lt, i->i^3));
           I := ideal(0_E); 
           for i from 1 to s do (
      	      h := new MutableList; 
      	      scan((#lt),i->h#i=randomIdeal'(randomIdeal(lt#i))); 
      	      I = I+sum(toList h); 
      	      );
           for i from 1 to 4 do (
      	      r := randomIdeal'(ideal(random(E^{0},E^{3:-1})));
      	      I = I+r;
      	      );
           rank source gens trim  (I+J) == t
           );
i18 : for i from 11 to 16 list secondStep(i)
o18 = {true, true, true, true, true, true}
o18 : List
\end{verbatim}
\end{footnotesize}

Now assume $n\ge 17$. Let $N$ be a third general codimension six subspace of $V$ and
let  $\N$ be the Grassmann variety of $3$-planes in $N$.

We have a short exact sequence of sheaves
$$
0\rightarrow I_{{ \sL}\cup{ \M}\cup{ \N},{\P {{}}(\wedge^3 V)}}(1)\rightarrow 
I_{{ \sL}\cup{\M},{\P {{}}(\wedge^3 V)}}(1)\rightarrow I_{\left({ \sL}\cup{ \M}\right)\cap { \N},\N}(1)\rightarrow 0. 
$$
To prove case (i), we will specialize the 4 points in
$\{r_i\}$ to lie on $\N$, $\lfloor\frac{6n-85}{9}\rfloor$ of the points in $\{p_i\}$ to lie on $\sL\cap \N$
and $\lfloor\frac{6n-85}{9}\rfloor$ of the points in $\{q_i\}$ to lie on $\M\cap \N$. This leaves in place
exactly four of the points in $\{p_i\}$ and four of the points in $\{q_i\}$.

Let $Y$ be the union of the tangent spaces at the points in
$\{p_i\}, \{q_i\}$ and  $\{r_i\}$. Then we obtain the following exact sequence: 
\[
0
\rightarrow  H^0(I_{\mathcal{K} \cup{ \N},{\P {{}}(\wedge^3 V)}}(1))
\rightarrow 
H^0(I_{\mathcal{K} ,{\P {{}}(\wedge^3 V)}}(1))
\rightarrow 
H^0(I_{\mathcal{K}\cap { \N}, \N}(1)). 
\]
where $\mathcal{K} = Y\cup { \sL}\cup{\M}$.
Note that we have the isomorphism: 
\[
H^0\left(I_{\mathcal{K} \cap { \N},\N}(1)\right) \simeq H^0\left( I_{\mathcal{K} \cap { \N},\P{{}}(\wedge^3 N)}(1)\right). 
\]
Thus the following inequality holds: 
\[
\dim H^0\left(I_{\mathcal{K} ,{\P {{}}(\wedge^3 V)}}(1) \right) \geq 
\dim H^0\left(I_{\mathcal{K} \cup{ \N},{\P {{}}(\wedge^3 V)}}(1)\right)+
\dim H^0\left( I_{\mathcal{K} \cap { \N},\P{{}}(\wedge^3 N)}(1)\right). 
\]
Since $Y\cap \N$ satisfies the conditions of Proposition~\ref{propa}, 
$H^0\left(I_{\mathcal{K} \cup{ \N},{\P {{}}(\wedge^3 V)}}(1)\right)$ 
has the expected dimension. 
By the induction hypothesis, 
$\dim H^0\left( I_{\mathcal{K}\cap { \N},\P{{}}(\wedge^3 N)}(1)\right)$ 
also has the expected value. Thus  
$\dim H^0\left(I_{\mathcal{K},{\P {{}}(\wedge^3 V)}}(1) \right)$ 
has the expected value, which proves (i).

 To prove case (ii) we make a similar specialization but substituting 
$\lfloor\frac{6n-85}{9}\rfloor$ with $\lceil\frac{6n-85}{9}\rceil$.
\end{proof}

Note that $f(n)-f(n-6)=3n^2-18n+35$.
In particular the system of hyperplanes in $\P{{}}(\wedge^3 V)$ which contain
$\sL$ has dimension $3n^2-18n+35$ for $n\ge 8$.
We want to keep $\lfloor\frac{6n-13}{9}\rfloor$ points outside $\sL$. Note that
$$\frac{3n^2-18n+35-(3n-5)(6n-13)/9}{18}=\frac{n^2}{18}-\frac{31n}{54}+\frac{125}{81}.$$

This leads us to the next proposition:
\begin{proposition}\label{propc}
Let $X=Gr(2,n)$ with $n\ge 9$. Let $V=\mathbb{K}^{n+1}$. Let $L$ be a general codimension six subspace of $V$ and
let $\sL$ be the Grassmann variety of $3$-planes in $L$. 
If 
\[ f_1(n):=\left\lfloor\frac{n^2}{18}-\frac{31n}{54}+\frac{125}{81}-\frac{n}{6}+2\right\rfloor
\]
and 
\[
f_2(n):=\left\lceil\frac{n^2}{18}-\frac{31n}{54}+\frac{125}{81}+\frac{n}{6}-1\right\rceil, 
\]
then 
\begin{itemize}
\item[(i)]
The system of hyperplanes in $\P{{}}(\wedge^3 V)$ which contain
$\sL$ and which contain the tangent spaces at $f_1(n)$
general points in $\sL$ and at $\lfloor\frac{6n-13}{9}\rfloor$
general points in $X$ has the expected dimension
$3n^2-18n+35-18f_1(n)-(3n-5)\lfloor\frac{6n-13}{9}\rfloor=O(n)$.

\item[(ii)] 
There are no hyperplanes in $\P{{}}(\wedge^3 V)$ which contain
$\sL$ and which contain the tangent spaces at $f_2(n)$
general points in $\sL$ and at $\lceil\frac{6n-13}{9}\rceil$
general points in $X$ .
\end{itemize}
\end{proposition}
\begin{proof}
We will let $\{p_i\}$ denote a set of $f_1(n)$ (resp. $f_2(n)$) general points on $\sL$ and let $\{q_i\}$ denote a set of $\lfloor\frac{6n-13}{9}\rfloor$ (resp. $\lceil\frac{6n-13}{9}\rceil$) general points on $X$. The proof is by a 6-step induction from $n-6$ to $n$. The initial cases
$n=9, 10, 11, 12, 13, 14$ can be checked directly as follows:
\begin{footnotesize}
\begin{verbatim}
i19 : thirdStep = method()
o19 = thirdStep
o19 : MethodFunction
i20 : thirdStep(ZZ) := n -> (
           f := floor(n^2/18-31*n/54+125/81-n/6+2);
           s := floor((6*n-13)/9); 
           t := binomial(n+1,3)-(3*n^2-18*n+35-18*f-(3*n-5)*s); 
           E := ZZ/32003[e_0..e_n,SkewCommutative=>true];
           L := ideal(e_6..e_n);
           Lk := L^3;
           I := ideal(0_E); 
           for i from 1 to f do (
                I = I+randomIdeal'(randomIdeal(L));
          	  );
           for i from 1 to s do (
                I = I+randomIdeal'(ideal(random(E^{1:0},E^{3:-1})));
      	      );
           rank source gens trim  (I+Lk) == t
           );
i21 : for i from 9 to 14 list thirdStep(i)
o21 = {true, true, true, true, true, true}
o21 : List
\end{verbatim}
\end{footnotesize}

Now assume that $n\ge 15$. Let $M$ be a second general codimension six subspace of $V$ and
let  $\M$ be the Grassmann variety of $3$-planes in $M$.

 We have the short exact sequence of sheaves
$$
0\rightarrow I_{{\sL}\cup{\M},{{\P {{}}(\wedge^3 V)}}}(1)\rightarrow 
I_{{\sL},{{\P {{}}(\wedge^3 V)}}}(1)\rightarrow I_{{\sL}\cap {\M}, \M}(1)\rightarrow 0. 
$$

To prove case (i)  we will specialize
$f_1(n)-\lfloor\frac{6n-49}{9}\rfloor$ of the
 points in $\{p_i\}$ to $\sL\cap \M$  and $\lfloor\frac{6n-49}{9}\rfloor$ of the
 points in $\{q_i\}$ to $\M$. This leaves in place exactly $\lfloor\frac{6n-49}{9}\rfloor$ of the
 points in $\{p_i\}$ and four of the points in $\{q_i\}$.

Let $Y$ be the union of the tangent spaces at the points in
$\{p_i\}$ and $\{q_i\}$. Then, as in the proof of the previous proposition, 
we obtain the following exact sequence: 
$$
0\rightarrow H^0 (I_{Y\cup{\sL}\cup{\M},{\P {{}}(\wedge^3 V)}}(1))\rightarrow 
H^0(I_{X\cup {\sL},{\P {{}}(\wedge^3 V)}}(1))
\rightarrow H^0(I_{(X\cup{\sL})\cap {\M},{\P {{}}(\wedge^3 M)}}(1)).
$$

By the induction hypothesis, the third element in the exact sequence has the expected dimension. 
Note that
$$
f_1(n)-\left\lfloor\frac{6n-49}{9}\right\rfloor\le f_1(n-6) 
$$
(in fact, the summand $-\frac{n}{6}$ was inserted in the definition of $f_1$ in order for this inequality to hold).
Now if we apply Proposition~\ref{propb} to the first element in the exact sequence, we can prove (i). 

To prove case (ii)  we will specialize
$f_2(n)-\lceil\frac{6n-49}{9}\rceil$ of the
 points in $\{p_i\}$ to $\sL\cap \M$ and $\lceil\frac{6n-49}{9}\rceil$ of the
 points in $q_i$ to $\M$. Since $$
f_2(n-6)\le f_2(n)-\left\lceil\frac{6n-49}{9}\right\rceil, 
$$
 we can apply Proposition~\ref{propb} to the first element in the exact sequence and we are done. 
\end{proof}

We now have the tools in place to prove the main theorem of this section (\thmref{asymp} of the Introduction).

\begin{theorem}
Let $n\ge 9$. Let $$s_1(n)=\left\lfloor\frac{n^2}{18}-\frac{2n}{27}+\frac{170}{81}\right\rfloor \ \ \mbox{and} 
\ \ s_2(n)=\left\lceil\frac{n^2}{18}+\frac{7n}{27}-\frac{73}{81}\right\rceil.$$
Then $\sigma_s(Gr(2,n))$ has the expected dimension whenever $s\le s_1(n)$
and whenever $s\ge s_2(n)$ (in this second case it fills the ambient space).
\end{theorem}

\begin{proof}
The proof is by a 6-step induction from $n-6$ to $n$. The cases
$n=9, 10, 11, 12, 13, 14$ can be checked directly and are well known (\cite{McG}). Let $n\ge 15$. Let $V=\mathbb{K}^{n+1}$. 
Let $L$ be a general codimension six subspace of $V$ and
let  $\sL$ be the Grassmann variety of $3$-planes in $L$. We will let $\{p_i\}$ denote a set of $s_1(n)$ (resp. $s_2(n)$) general points on $Gr(2,n)$.
Note that $$
s_1(n)=f_1(n)+\left\lfloor\frac{6n-13}{9}\right\rfloor \ \ \mbox{and} 
\ \ s_2(n)=f_2(n)+\left\lceil\frac{6n-13}{9}\right\rceil.
$$

Consider the following short exact sequence of vector spaces:
$$
0\rightarrow H^0(I_{{\sL},{{\P {{}}(\wedge^3 V)}}}(1))\rightarrow 
\wedge^3 V\rightarrow \wedge^3 L\rightarrow 0.
$$
To prove that $\sigma_s(Gr(2,n))$ has the expected dimension whenever $s\le s_1(n)$, we specialize $f_1(n)$ of the points in $\{p_i\}$ to lie on $\sL$
and we keep $\left\lfloor\frac{6n-13}{9}\right\rfloor$ points in their place.
Let $Y$ be the union of the tangent spaces to $X$ at the points in $\{p_i\}$. 
Then we obtain the following exact sequence:  
$$
0\rightarrow  H^0(I_{Y\cup{\sL},{\P {{}}(\wedge^3 V)}}(1))\rightarrow 
H^0(I_{Y,{\P {{}}(\wedge^3 V)}}(1))
\rightarrow 
H^0(I_{Y\cap {\sL},{\P {{}}(\wedge^3 L)}}(1)). 
$$
By the induction hypothesis, $H^0(I_{Y\cap {\sL},{\P {{}}(\wedge^3 L)}}(1))$ has the expected dimension. 
Note that the following inequality holds: 
$$
s_1(n)-s_1(n-6)\le \left\lfloor\frac{6n-13}{9}\right\rfloor. 
$$
It follows from \propref{propc}  that $H^0(I_{Y\cup{\sL},{\P {{}}(\wedge^3 V)}}(1))$ also has 
the expected dimension. Thus we have proved that $\sigma_s(Gr(2,n))$ has the expected dimension whenever $s\le s_1(n)$. 
 
Since the following inequality holds: 
$$
\left\lceil\frac{6n-13}{9}\right\rceil\le s_2(n)-s_2(n-6),
$$
the proof that $\sigma_s(Gr(2,n))$ has the expected dimension whenever $s\ge s_2(n)$ can be shown in the same way by 
 specializing $f_2(n)$ of the points in $\{p_i\}$ to lie on $\sL$ and keeping $\left\lceil\frac{6n-13}{9}\right\rceil$ points in their place.
\end{proof}

\section{The defective cases}
In \conjref{BDdG} there is a list of four defective secant varieties of Grassmannians. All four of the defective cases  are described in \cite{CGG1}. We make here some further comments.

A geometric explanation of the defectivity of $X=Gr(3,7)$ is the following, inspired by \cite{CC2}. As in \cite{CGG1},
given three points $P_1, P_2, P_3$ in $X$, there is a basis $e_0,\ldots, e_7$ such that
the three points correspond to $P_1=\langle e_0,e_1,e_2,e_3 \rangle, P_2= \langle e_4,e_5,e_6,e_7 \rangle$, $P_3= \langle e_0+e_4,e_1+e_5,e_2+e_6,e_3+e_7 \rangle$.  Using the matrix

\[\left[\begin{array}{cccc|cccc}
1&&&&t\\
&1&&&&t\\
&&1&&&&t\\
&&&1&&&&t\\
\end{array}\right],\]
we see that there is a rational normal curve embedded with $\sO(4)$ which passes through the 3 points and is contained in $X$.
The existence of this curve, $C_4$, forces each of the tangent spaces $T_{P_i}X$ to have the line $T_{P_i}C_4$  in common with the
$\P 4$ spanned by $C_4$.  This leads to the following inequalities:
\[
\dim \langle T_{P_1}G, T_{P_2} G, T_{P_3} G \rangle \leq 4+3(\dim Gr(3,7)-1)=4 +3\cdot 15 = 49 < 50. 
\]
By Terracini's lemma, this proves the defectivity of $Gr(3, 7)$.
The defectivity of $\sigma_4(Gr(3, 7))$ follows as a direct consequence of the defectivity
 of $\sigma_3(Gr(3, 7))$. 
 
A geometric explanation of the defectivity of  $X=Gr(2,8)$ is similar.
 For any $4$ general points in
$X$ we find a Veronese surface embedded with $\sO(3)$ which passes through the 4 points and is contained in $X$.
Let $P_1, P_2, P_3, P_4$ be the four points.
We may assume that (see \cite{CGG1})
$P_1=\langle e_0,e_1,e_2 \rangle, P_2=\langle e_3,e_4,e_5 \rangle,
P_3=\langle e_6,e_7,e_8 \rangle, P_4= \langle e_0+e_3+e_6, e_1+e_4+e_7,e_2+e_5+e_8\rangle$.
If $(s,t,u)$ are projective coordinates of $\P2$ we use the matrix

\[\left[\begin{array}{ccc|ccc|ccc}
s&&&t&&&u\\
&s&&&t&&&u\\
&&s&&&t&&&u\\
\end{array}\right]\]
to realize the Veronese surface passing through the 4 points.
It follows that the span of the tangent spaces has dimension at
 most ${5\choose{2}}-1+4\cdot(18-2)=73$ while the expected dimension
is $4\cdot 18+3=75$. By Terracini's lemma, this proves the defectivity of $\sigma_4(Gr(2, 8))$.

The geometric argument for the defectivity of $X=Gr(2,6)$ is more subtle. It
can be proved by the following argument which also helps to find
the (set theoretical) equations of the secant varieties of $Gr(2,6)$.
Given $\omega\in\wedge^3{\mathbb C}^7$, there is a well defined contraction operator
$$\phi_{\omega}\colon\wedge^2{\mathbb C}^7\to\wedge^5{\mathbb C}^7. $$
Let $e_1,\ldots ,e_7$ be a basis of ${\mathbb C}^7$. If $\omega=e_1\wedge e_2\wedge e_3$
then  $\phi_{\omega}(e_i\wedge e_j)$ is nonzero if and only if $i,j\ge 4$. It follows
that $\textrm{rank}(\phi_{\omega})={4\choose 2}=6$. Since  $\textrm{rank}(\phi_{\omega})=\mathrm{rank}(\phi_{g\omega})$
for every $g\in SL({\mathbb C}^7)$, we get
that  $\textrm{rank}(\phi_{\omega})=6$ if $\omega\in Gr(2,6)$.

If $\omega=\sum_{i=1}^k \omega_i$ with $\omega_i$ decomposable (i.e. $\omega_i\in Gr(2,6)$) then it follows that
$\textrm{rank}(\phi_{\omega})=\textrm{rank}(\sum_{i=1}^k\phi_{\omega_i})\le \sum_{i=1}^k\textrm{rank}(\phi_{\omega_i})\le 6k$.
Hence if $\omega\in \sigma_k(Gr(2,6))$ then by semicontinuity we have $\textrm{rank}(\phi_{\omega}) \le 6k$.

Consider $\omega=e_1\wedge e_3\wedge e_5+e_1\wedge e_4\wedge e_7+e_1\wedge e_2\wedge e_6+ e_2\wedge e_3\wedge e_4
+e_5\wedge e_6\wedge e_7.$
We can represent $\omega$ via the diagram

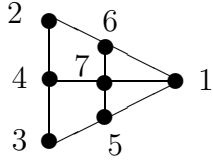
\begin{figure}[h]
\begin{center}
\setlength{\unitlength}{3947sp}%
\begingroup\makeatletter\ifx\SetFigFont\undefined%
\gdef\SetFigFont#1#2#3#4#5{%
  \reset@font\fontsize{#1}{#2pt}%
  \fontfamily{#3}\fontseries{#4}\fontshape{#5}%
  \selectfont}%
\fi\endgroup%
\begin{picture}(735,1280)(4000,-800)
\thinlines
\put(4255,-676){\circle*{100}}
\put(4255, 74){\circle*{100}}
\put(4604,-526){\circle*{100}}
\put(4259,-286){\circle*{100}}
\put(5050,-301){\circle*{100}}
\put(4604,-316){\circle*{100}}
\put(4610,-91){\circle*{100}}
\put(5080,-301){\line(-2,-1){810}}
\put(5080,-301){\line(-2, 1){810}}
\put(4255,-676){\line( 0, 1){750}}
\put(5020,-301){\line(-1, 0){795}}
\put(4606,-91){\line( 0,-1){450}}
\put(5191,-361){\makebox(0,0)[lb]{\smash{\SetFigFont{12}{14.4}{\rmdefault}{\mddefault}{\updefault}1}}}
\put(4021,-346){\makebox(0,0)[lb]{\smash{\SetFigFont{12}{14.4}{\rmdefault}{\mddefault}{\updefault}4}}}
\put(4411,-271){\makebox(0,0)[lb]{\smash{\SetFigFont{12}{14.4}{\rmdefault}{\mddefault}{\updefault}7}}}
\put(4591, 14){\makebox(0,0)[lb]{\smash{\SetFigFont{12}{14.4}{\rmdefault}{\mddefault}{\updefault}6}}}
\put(3991, 59){\makebox(0,0)[lb]{\smash{\SetFigFont{12}{14.4}{\rmdefault}{\mddefault}{\updefault}2}}}
\put(4021,-736){\makebox(0,0)[lb]{\smash{\SetFigFont{12}{14.4}{\rmdefault}{\mddefault}{\updefault}3}}}
\put(4621,-766){\makebox(0,0)[lb]{\smash{\SetFigFont{12}{14.4}{\rmdefault}{\mddefault}{\updefault}5}}}
\end{picture}
\end{center}
\end{figure}

\noindent An explicit computation shows that $rank(\phi_{\omega})=21$. It follows that 
$\sigma_3(Gr(2,6))$ cannot fill the ambient space, hence $Gr(2,6)$ is defective.

\begin{theorem}
Let $\omega\in\wedge^3{\mathbb C}^7$. Consider the contraction operator $\phi_{\omega}\colon\wedge^2{\mathbb C}^7\to\wedge^5{\mathbb C}^7. $
The equation of $\sigma_3(Gr(2,6))$ is given by an $SL(7)$-invariant polynomial $P_7$
of degree seven such that
$$\det(\phi_{\omega})=2\left[P_7(\omega)\right]^3.$$
\end{theorem}

\begin{proof}
The morphism $\phi_{\omega}$ drops rank by three when $\omega$ belongs
to the hypersurface $\sigma_3(Gr(2,6))$. Hence the linear embedding of
$\wedge^3{\mathbb C}^7$ in $\mathrm{Hom}(\wedge^2{\mathbb C}^7,\wedge^5{\mathbb C}^7) $
(given by $\omega\mapsto \phi_{\omega}$)
meets the determinantal hypersurface with multiplicity three. By direct computation
on $\omega=a_{135}e_1\wedge e_3\wedge e_5+a_{147}e_1\wedge e_4\wedge e_7+a_{126}e_1\wedge e_2\wedge e_6+ a_{234}e_2\wedge e_3\wedge e_4
+a_{567}e_5\wedge e_6\wedge e_7$, we see that
$$\det\phi_{\omega}=-2(a_{234}^2a_{567}^2a_{135}a_{147}a_{126})^3.$$ Hence we can arrange the scalar multiples
in order that $P$ is defined over the rational numbers and the equation $\det(\phi_{\omega})=2\left[P_7(\omega)\right]^3$
holds. 
\end{proof}

\begin{remark} The equation $v\wedge v'\wedge\omega=v'\wedge v\wedge\omega$ for
$v, v'\in\wedge^2{\mathbb C}^7$ shows that $\phi_{\omega}$ is symmetric.
A natural symmetric operator such that its determinant is a cube appears already in
\cite{Ot}, where the coefficient $2$ appears at the same place. 
 The coefficient
$2$ is needed if we want the invariant $P_7$ to be defined over the rational numbers.

The graphical notation found in the above diagram comes from the original paper of Schouten \cite{Sch}.
Indeed, the case $Gr(2,6)$ is in principle well known because $SL(7)$ has only finitely many orbits
on $\P{{}} (\wedge^3{\mathbb C}^7)$. This classification was computed in 1931 by Schouten, \cite{Sch},
correcting previous work of Reichel, who missed the orbit of dimension $20$.
 He found 
 all of the 9 orbits for this 
action together with their dimensions.

G.B. Gurevich in his textbook \cite{Gu} gave equations for these orbits but  from his 
description it is not easy to find the order among the orbits.
 In 
fact the obvious order relation (Bruhat order) among the orbits, such that $ 
O_1 
\leq O_2 
$ if the closure of $ O_2 $ contains $ O_1 $, is {\it not a total order}, 
and indeed this is the first case among Grassmannians where this phenomenon occurs.
We take the opportunity to show in the following table 
the order relation among the 9 orbits, computed by Elisabetta Ardito in her laurea thesis, defended in 
L'Aquila in 1997 under the supervision of the second author.
It is a distributive lattice.
 We have added to each orbit the value of
$\mathrm{rank}(\phi_{\omega})$ together with some geometrical information. Each of the values can be computed easily
on a representative of each orbit. The dimensions can be computed by considering the rank of the derivative
of the action of $SL(7)$ on $\wedge^3{\mathbb C}^7$. 
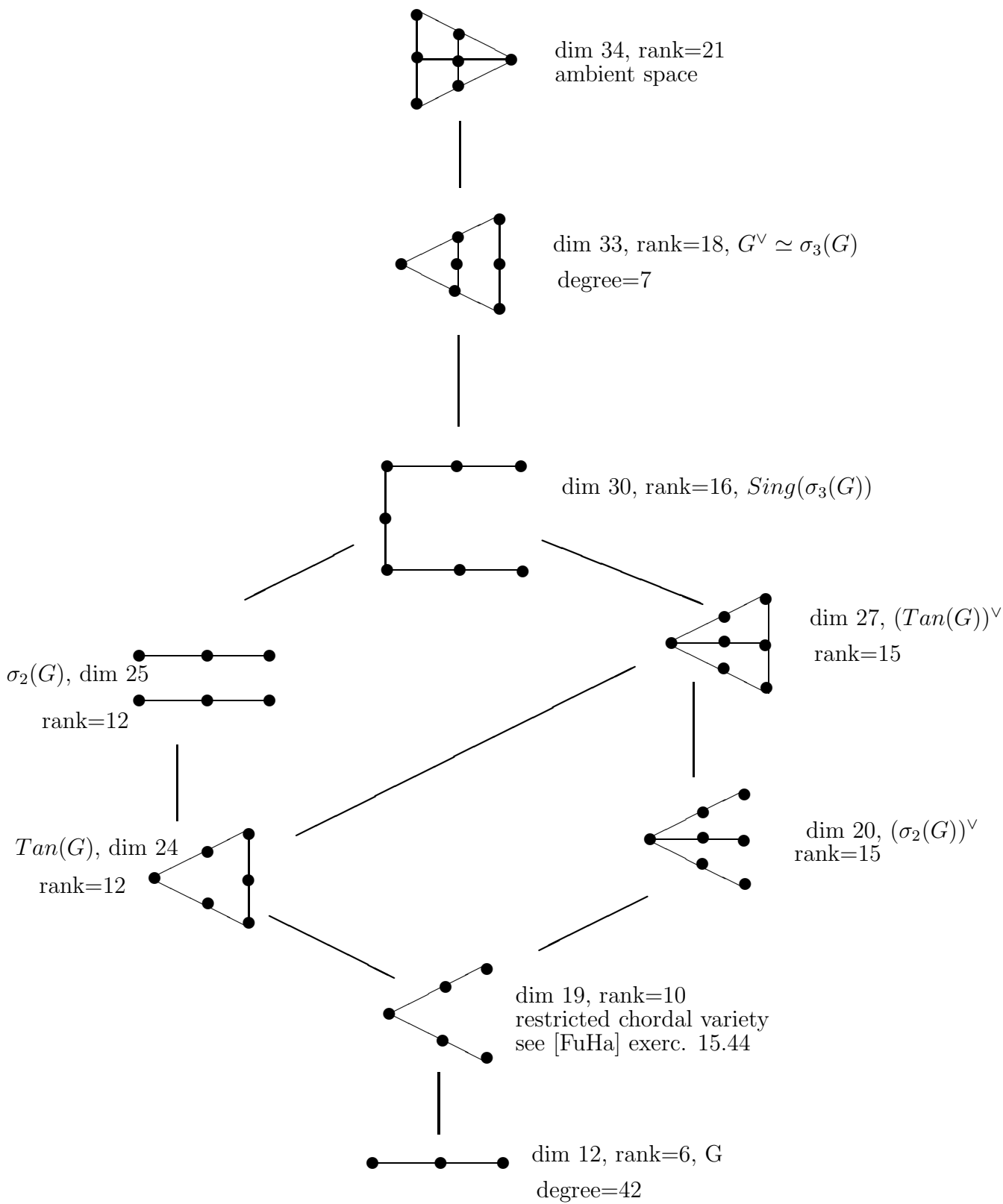
\begin{figure}
\begin{center}
\setlength{\unitlength}{3947sp}%
\begingroup\makeatletter\ifx\SetFigFont\undefined%
\gdef\SetFigFont#1#2#3#4#5{%
  \reset@font\fontsize{#1}{#2pt}%
  \fontfamily{#3}\fontseries{#4}\fontshape{#5}%
  \selectfont}%
\fi\endgroup%
\begin{picture}(6360,10238)(1401,-10044)
\thinlines
\put(4126,-2026){\circle*{100}}
\put(4951,-1651){\circle*{100}}
\put(4951,-2401){\circle*{100}}
\put(4576,-2251){\circle*{100}}
\put(4602,-1801){\circle*{100}}
\put(4951,-2026){\circle*{100}}
\put(4126,-2026){\line( 2, 1){810}}
\put(4126,-2026){\line( 2,-1){810}}
\put(4951,-1651){\line( 0,-1){750}}
\put(4606,-1816){\line( 0,-1){465}}
\put(4591,-2026){\circle*{100}}
\put(4255,-676){\circle*{100}}
\put(4255, 74){\circle*{100}}
\put(4604,-526){\circle*{100}}
\put(4259,-286){\circle*{100}}
\put(5050,-301){\circle*{100}}
\put(4604,-316){\circle*{100}}
\put(4610,-91){\circle*{100}}
\put(5080,-301){\line(-2,-1){810}}
\put(5080,-301){\line(-2, 1){810}}
\put(4255,-676){\line( 0, 1){750}}
\put(5020,-301){\line(-1, 0){795}}
\put(4606,-91){\line( 0,-1){450}}
\thicklines
\put(4606,-2626){\line( 0,-1){765}}
\put(4621,-826){\line( 0,-1){555}}
\thinlines
\put(3882,-9601){\circle*{100}}
\put(4981,-9601){\circle*{100}}
\put(4456,-9601){\circle*{100}}
\put(3856,-9601){\line( 1, 0){1125}}
\put(1917,-5701){\circle*{100}}
\put(3016,-5701){\circle*{100}}
\put(2491,-5701){\circle*{100}}
\put(1891,-5701){\line( 1, 0){1125}}
\put(1917,-5326){\circle*{100}}
\put(3016,-5326){\circle*{100}}
\put(2491,-5326){\circle*{100}}
\put(1891,-5326){\line( 1, 0){1125}}
\put(2842,-6826){\circle*{100}}
\put(2842,-7576){\circle*{100}}
\put(2493,-6976){\circle*{100}}
\put(2838,-7216){\circle*{100}}
\put(2497,-7411){\circle*{100}}
\put(2047,-7201){\circle*{100}}
\put(2017,-7201){\line( 2, 1){810}}
\put(2017,-7201){\line( 2,-1){810}}
\put(2842,-6826){\line( 0,-1){750}}
\put(4021,-8341){\circle*{100}}
\put(4846,-7966){\circle*{100}}
\put(4846,-8716){\circle*{100}}
\put(4471,-8566){\circle*{100}}
\put(4497,-8116){\circle*{100}}
\put(4021,-8341){\line( 2, 1){810}}
\put(4021,-8341){\line( 2,-1){810}}
\put(7012,-6496){\circle*{100}}
\put(6663,-6646){\circle*{100}}
\put(7008,-6886){\circle*{100}}
\put(6217,-6871){\circle*{100}}
\put(6663,-6856){\circle*{100}}
\put(6657,-7081){\circle*{100}}
\put(7017,-7246){\circle*{100}}
\put(6187,-6871){\line( 2, 1){810}}
\put(6187,-6871){\line( 2,-1){810}}
\put(6247,-6871){\line( 1, 0){795}}
\put(7192,-4846){\circle*{100}}
\put(6843,-4996){\circle*{100}}
\put(7188,-5236){\circle*{100}}
\put(6397,-5221){\circle*{100}}
\put(6843,-5206){\circle*{100}}
\put(6837,-5431){\circle*{100}}
\put(7197,-5596){\circle*{100}}
\put(6367,-5221){\line( 2, 1){810}}
\put(6367,-5221){\line( 2,-1){810}}
\put(6427,-5221){\line( 1, 0){795}}
\put(7216,-4816){\line( 0,-1){825}}
\put(4006,-4605){\circle*{100}}
\put(5146,-4616){\circle*{100}}
\put(4617,-4601){\circle*{100}}
\put(4021,-4601){\line( 1, 0){1125}}
\put(4006,-3731){\circle*{100}}
\put(5131,-3731){\circle*{100}}
\put(4591,-3731){\circle*{100}}
\put(4006,-3731){\line( 1, 0){1125}}
\put(3991,-4171){\circle*{100}}
\put(3991,-3766){\line( 0,-1){825}}
\thicklines
\put(4441,-8836){\line( 0,-1){525}}
\put(3721,-4396){\line(-2,-1){912}}
\put(2236,-6076){\line( 0,-1){645}}
\put(6586,-5551){\line( 0,-1){795}}
\put(3016,-7516){\line( 2,-1){1050}}
\put(3241,-6841){\line( 2, 1){2850}}
\put(6196,-7351){\line(-2,-1){912}}
\put(5311,-4351){\line( 5,-2){1350}}
\put(5221,-9586){\makebox(0,0)[lb]{\smash{\SetFigFont{12}{14.4}{\rmdefault}{\mddefault}{\updefault}dim 12, rank=6, G}}}
\put(5261,-9886){\makebox(0,0)[lb]{\smash{\SetFigFont{12}{14.4}{\rmdefault}{\mddefault}{\updefault}degree=42}}}
\put(5086,-8251){\makebox(0,0)[lb]{\smash{\SetFigFont{12}{14.4}{\rmdefault}{\mddefault}{\updefault}dim 19, rank=10}}}
\put(5086,-8451){\makebox(0,0)[lb]{\smash{\SetFigFont{12}{14.4}{\rmdefault}{\mddefault}{\updefault}restricted chordal 
variety}}}
\put(5086,-8651){\makebox(0,0)[lb]{\smash{\SetFigFont{12}{14.4}{\rmdefault}{\mddefault}{\updefault}see \cite{FuHa} exerc. 
15.44}}}
\put(7561,-5071){\makebox(0,0)[lb]{\smash{\SetFigFont{12}{14.4}{\rmdefault}{\mddefault}{\updefault}dim 27, 
$(Tan(G))^{\vee}$}}}
\put(7600,-5371){\makebox(0,0)[lb]{\smash{\SetFigFont{12}{14.4}{\rmdefault}{\mddefault}{\updefault}rank=15}}}
\put(5401,-1921){\makebox(0,0)[lb]{\smash{\SetFigFont{12}{14.4}{\rmdefault}{\mddefault}{\updefault}dim 33, rank=18, 
$G^{\vee}\simeq \sigma_3(G)$}}}
\put(5440,-2221){\makebox(0,0)[lb]{\smash{\SetFigFont{12}{14.4}{\rmdefault}{\mddefault}{\updefault}degree=7}}}
\put(5416,-285){\makebox(0,0)[lb]{\smash{\SetFigFont{12}{14.4}{\rmdefault}{\mddefault}{\updefault}dim 34, rank=21}}}
\put(5416,-485){\makebox(0,0)[lb]{\smash{\SetFigFont{12}{14.4}{\rmdefault}{\mddefault}{\updefault}ambient space}}}
\put(876,-7005){\makebox(0,0)[lb]{\smash{\SetFigFont{12}{14.4}{\rmdefault}{\mddefault}{\updefault}$Tan(G)$, dim 24}}}
\put(1076,-7325){\makebox(0,0)[lb]{\smash{\SetFigFont{12}{14.4}{\rmdefault}{\mddefault}{\updefault}rank=12}}}
\put(7531,-6856){\makebox(0,0)[lb]{\smash{\SetFigFont{12}{14.4}{\rmdefault}{\mddefault}{\updefault}dim 20, 
$(\sigma_2(G))^{\vee}$}}}
\put(7431,-7056){\makebox(0,0)[lb]{\smash{\SetFigFont{12}{14.4}{\rmdefault}{\mddefault}{\updefault}rank=15}}}
\put(801,-5536){\makebox(0,0)[lb]{\smash{\SetFigFont{12}{14.4}{\rmdefault}{\mddefault}{\updefault}$\sigma_2(G)$, dim 25}}}
\put(1100,-5936){\makebox(0,0)[lb]{\smash{\SetFigFont{12}{14.4}{\rmdefault}{\mddefault}{\updefault}rank=12}}}
\put(5476,-3961){\makebox(0,0)[lb]{\smash{\SetFigFont{12}{14.4}{\rmdefault}{\mddefault}{\updefault}dim 30, rank=16,
$Sing(\sigma_3(G))$}}}
\end{picture}
\caption{Orbits for $SL(7)$-action on $\wedge^3{\mathbb C}^7$}
\label{figure}
\end{center}
\end{figure}
\end{remark}
It follows from this description the following theorem:

\begin{theorem}
For $\omega\in \wedge^3{\mathbb C}^7$ the following holds:
\begin{itemize}
\item[$\mathrm{(i)}$]  $\omega\in Gr(2,6)$ if and only if  $\mathrm{rank}(\phi_{\omega})\le 6$.
\item[$\mathrm{(ii)}$]  $\omega\in \sigma_2(Gr(2,6))$ if and only if  $\mathrm{rank}(\phi_{\omega})\le 12$. Hence
the $13\times 13$ minors of $\phi_{\omega}$ give set theoretic equations of
$\sigma_2(Gr(2,6))$.
\item[$\mathrm{(iii)}$] $\omega\in \sigma_3(Gr(2,6))$ if and only if  $\mathrm{rank}(\phi_{\omega})\le 18$.
\end{itemize}
\end{theorem}

In particular the table on the following page shows the possible degenerations of elements in $ 
\wedge ^3 
{\mathbb C}^7$.
There are two degenerations which are not obvious:

\begin{itemize}
\item 
The degeneration of P$_{25}$ in P$_{24}$ (the subscript means the dimension) resulting from
\[
\begin{array}{l}
 \lim_{t \rightarrow 0}\frac{1}{t}[e_1 \wedge e_2 \wedge e_3 - (e_1+te_4)\wedge 
(e_2+te_5) \wedge
(e_3+te_6)] \\
= -(e_1 \wedge e_2 \wedge e_6+e_1 \wedge e_5 \wedge e_3+e_4 \wedge e_2 \wedge e_3)
\end{array}
\]
\item 
The degeneration of P$_{30}$ in P$_{27}$ resulting from
\[
\begin{array}{ll}
 \lim_{t \rightarrow 0}\frac{1}{t}[(e_1 \wedge e_2 \wedge e_3+(e_3 \wedge e_7 \wedge(e_3+te_6)) &  \\
-(e_3+te_6)\wedge(e_1+te_4)\wedge(e_2+te_5)] & \\
= e_3 \wedge e_7 \wedge e_6 - e_3 \wedge e_1 \wedge e_5 - e_3 \wedge e_4 \wedge e_2 - 
e_6 \wedge e_1 \wedge e_2 & 
\end{array}
\]
\end{itemize}

All the other degenerations are somewhat clear by considering the shape of the diagrams.
Vinberg's school \cite{VE} extended Schouten and Gurevich's classification to higher dimension, but 
the Bruhat order  of the orbits has not yet been explicitly written.

Notice that  the hypersurface $\sigma_3(Gr(2,6))$   is isomorphic to the dual variety
of  $Gr(2,6)$ (which has degree $7$, see \cite{La}).
It is called $C_8$ in the notation of \cite{Gu}, pg. 393. 

It can be computed that $P_7$ is a polynomial with 10,680 terms.
The ideal of the secant variety $\sigma_2(Gr(2,6))$
is generated by 28 cubics which correspond to the ideal
$\Gamma^{3,1^6}V\subset S^3(\wedge^3{\mathbb C}^7)$, which is the covariant $C_4$
according to \cite{Gu}, pg. 393.
It is interesting to check that $\mathrm{Sing}(\sigma_3(Gr(2,6)))$ is the orbit of dimension $30$,
while $\mathrm{Sing}(\sigma_2(Gr(2,6)))=Gr(2,6)$.




\end{document}